\documentstyle{amsart}
 \newtheorem{thm}{Theorem}[section]
 
 \newtheorem{cor}[thm]{Corollary}

 \theoremstyle{definition}
 \newtheorem{rem}[thm]{Remark}
 
 \newtheorem{defn}[thm]{Definition}
 \newcommand{\key}{\bibitem}
 \newcommand{\red}{{\rm red}}

\begin{document}
 \title[Approximation by smooth curves near the tangent
 cone]{Approximation by smooth curves \\ near the tangent cone}
 \author{Anvar R. Mavlyutov}
 \address {Department of Mathematics, Indiana University, Bloomington, IN 47405, USA.}
 \email{amavlyut@@indiana.edu}

 \subjclass{Primary: 14}

 \begin{abstract}
 We show that through a  point of an affine variety
 there  always exists a smooth plane curve inside the ambient affine space,
 which has the multiplicity of intersection with the variety at least 3.
 This result has an application to the study of affine schemes.
 \end{abstract}

 \maketitle

 This note appeared as a result of studying certain affine schemes.
 Namely, I.~R. Shafarevich in \cite{sh}
 studied the question: which irreducible components of the reduced scheme
 $C_n^{\red}$
 of associative commutative multiplications on a $n$-dimensional vector space
 consist of those multiplications
 that represent nilpotent degree 3 algebras?
 While this problem was partially solved in \cite{sh} and \cite{am}, the methods
 used there
 do not suffice to give a complete answer to the question.
 The basic idea, due to Shafarevich, was to compare the tangent spaces to the
 non-reduced
 scheme $C_n$ and the scheme $A_n$ of nilpotent degree 3 commutative multiplications
 on a $n$-dimensional vector space.
 In certain cases, the  difference of the tangent vectors
 was explained by the presence of nilpotents in the structure sheaf
 of the scheme $C_n$  which can be
 eliminated by  embedding $C_n$ into the scheme of multiplications on a
 $(n+1)$-dimensional
 vector space representing
 associative algebras with a unit.
 Another method, used in \cite{am}, is to apply ``obstructions to deformations'' of
 algebras.
 However, this approach required the variety $C_n^\red$  to be  smooth at a point of
 $A_n^\red$ (a priori, we don't know if this is true)
 in order to deduce a definite answer to Shafarevich's question.
 The result of this paper shows that we can still use the second order obstructions
 to deformations of algebras, even if a whole irreducible component
 of $A_n^\red$ is in the singular locus of $C_n^{\red}$.

 Our plan is  to introduce the multiplicity of intersection
 of a smooth curve with a variety at a point  in Section~\ref{s:pc}, and then 
 show that an affine variety
 can be nicely ``approximated'' by smooth plane curves inside the ambient affine
 space.
 Section~\ref{s:ap} discusses an application of Theorem~\ref{t:m} to the studying
 of  affine schemes of algebras.

 \section{Plane curves approximating varieties near the tangent cone.}\label{s:pc}

 In this section, we generalize the notion of multiplicity of intersection
 used in \cite[Chapter~3,~\S4, Definition~1]{clo}.  Proposition~3 in
 \cite[Chapter~9,~\S6]{clo},
 has a criterion which says that a line lies in the tangent space
 of the variety iff the line meets the variety with multiplicity at least 2.
 Here, we show that through a singular point of an affine variety
 we can always draw a smooth plane curve inside the ambient affine space,
 which has the multiplicity of intersection with the variety at least 3.
 This basic result does not seem to appear anywhere in
 the literature.

 Let $K$ be an arbitrary field.
 An algebraic curve $C\subset{\Bbb A}^{n}$ can be
 parameterized at its smooth point $p$ by Taylor's series
 $$G(t)=p+t\cdot v_1+t^2\cdot v_2+\cdots,$$
 where $t$ is a local parameter.
 Using this parameterization, we can define multiplicity of intersection
 of an affine variety with the curve $C$ at the point $p$.

 \begin{defn}
 Let $m$ be a nonnegative integer and $X\subset{\Bbb A}^{n}$
 an affine variety, given by
 $\text{\bf I}(X)\subset K[x_1,\ldots,x_n]$. Suppose that we have
 a curve $C$, smooth at a point $p$ and with the parameterization as above.
 Then $C$ {\bf meets} $X$ {\bf with multiplicity} $m$  at $p$
 if $t=0$ is a zero of multiplicity at least $m$ of the series
 $f\circ G(t)$  for all $f\in\text{\bf I}(X)$, and for some $f$, of multiplicity
 exactly $m$. We denote this multiplicity of intersection by $I(p,C,X)$.
 \end{defn}

 \begin{rem}
 This definition does not depend on the choice of the parameterization.
 To check this one can use the following criterion:
 $t=0$ is a zero of multiplicity $m$ of $h(t)\in K[[t]]$ if and only if
 $h(0)=h'(0)=\cdots=h^{(m-1)}(0)=0$ but $h^{(m)}(0)\ne 0$.
 \end{rem}

 Also, we denote by $T_{p}X$ ($TC_{p}X$)
 the tangent space (cone) of a variety $X$ at a point $p$.
 The next result contains a necessary condition for a line to be in the tangent
 cone.
 Unfortunately, it is not a sufficient one.

 \begin{thm}\label{t:m}
 Let $X\subset{\Bbb A}^{n}$ be an affine variety and $p\in X$ a singular
 point. Then for each $v\in TC_{p}X$ there exists a smooth
 curve $C\subset{\Bbb A}^{n}$ through $p$ such that
 $v\in T_{p}C$ and $I(p,C,X)\geq3$.
 \end{thm}

 \begin{pf} Let  $v\in TC_{p}X$.
 We can choose affine coordinates $x_{1},\ldots,x_{n}$, so that
 $p=(0,\ldots,0)$. Let $C$ be a curve through
 $p$ parameterized
 by $G(t)=t\cdot v+t^{2}\cdot\gamma$, where
 $v\in TC_{p}X$,
 $\gamma=(\gamma_1,\ldots,\gamma_{n})$.
 It suffices  to show that there is $\gamma$, which is  zero or not a scalar
 multiple
 of $v$, such that $t=0$ is a zero
 of multiplicity $\geq3$ of $f\circ G(t)$  for all $f\in\text{\bf I}(X)$.
 In this case, $C$ will be a smooth curve with the local parameter $t$
 at $p$.

 For each $f\in\text{\bf I}(X)$, denote by
 $l_{f}(x)=l_{1}x_{1}+\cdots+l_{n}x_{n}$
  and
 $q_{f}(x)$ the  linear and quadratic parts of $f$.
 According to this notation we define
 $$W=\{(l_{1},\ldots,l_{n},q_{f}(v)):\,f\in\text{\bf I}(X)\}\subset
 k^{n+1},$$ which is a  subspace  since $\text{\bf I}(X)$ is. For all
 $f\in\text{\bf I}(X)$,    we have
 $$f(t\cdot v+t^{2}\cdot\gamma)=l_{f}(t\cdot v)+l_{f}(t^{2}\cdot\gamma)+
 q_{f}(t\cdot v)+\text{ terms of degree }\geq3\text{ in }t.$$
 Since $v\in T_{p}X$, we get $l_{f}(v)=0$ for any $f\in\text{\bf I}(X)$.
 Hence,
 $$f(t\cdot v+t^{2}\cdot\gamma)\equiv
 t^{2}(l_{1}\gamma_{1}+\cdots+l_{n}\gamma_{n}+q_{f}(v))$$
 modulo terms of degree $\geq3$ in $t$. So we need to resolve
 equations $a_{1}\gamma_{1}+\cdots+a_{n}\gamma_{n}+a_{n+1}=0$ in variables
 $\gamma_{1},\ldots,\gamma_{n}$ for all possible
 $(a_{1},\ldots,a_{n+1})\in W$.

 Consider the scalar product $k^{n+1}\times k^{n+1}\rightarrow k$,
 which sends $a\times b$ to $a\cdot b:=
 \sum^{n+1}_{i=1}a_{i}b_{i}$.
 Using this scalar product, we will find the vector
 $(\gamma_{1},\ldots,\gamma_{n},1)\in W^{\perp}:=
 \{b\in k^{n+1}:\,b\cdot a=0\text{ for all } a\in W\}$.
 If $W\subseteq k^{n}\times\{0\}$, then we can put $\gamma=(0,\ldots,0)$
 to obtain the required curve.
 Otherwise, consider the  projection
 $\pi:\,k^{n+1}\rightarrow k^{n}\times\{0\}$, which assigns
 $(a_{1},\ldots,a_{n},0)$
  to
 $(a_{1},\ldots,a_{n+1})$.
 And, let $\tilde{\pi}:\,W\rightarrow k^{n}\times\{0\}$ be induced by $\pi$.
 The projection  $\tilde{\pi}$ is injective, because $v\in TC_{p}X$ and
 $l_{f}=0$ imply $q_{f}(v)=0$.
 If we denote $\widetilde{W}=\tilde{\pi}(W)$, then
 $\dim W^{\perp}=\dim\widetilde{W}^{\perp}$, by injectivity of $\tilde{\pi}$.
 Now, if $W^{\perp}$ is of the form $W_{1}\times\{0\}$,
 then by construction $\widetilde{W}^{\perp}=W_{1}\times k$,
 contradicting with the equality of the dimensions.
 Therefore, there exists a vector
 $(\gamma_{1},\ldots,\gamma_{n},1)\in W^{\perp}$,
 i.e., $l_{f}(\gamma)+q_{f}(v)=0$ for all $f\in I(X)$.
 We claim that this $\gamma$ is linearly independent of $v$.
 Indeed, since $W$ is not included into $k^{n}\times\{0\}$, we get
 $q_{f}(v)\ne0$
  for
 some $f\in\text{\bf I}(X)$, which implies $l_{f}(\gamma)\ne0$.
 On the other hand, since $v\in T_{p}X$, we get $l_{f}(v)=0$. Thus, we have found
 the desired $\gamma$.
 \end{pf}

 \begin{rem}
 This result cannot be improved in the following sense. For the curve $X$,
 given by the equation $x^{2}=y^{3}$,
 in the affine plane, and $p=(0,0)$
 there is no smooth curve $C$ such that $I(p,C,X)\geq4$.
 \end{rem}

 \begin{rem}
 In the case of analytic varieties the theorem is also valid.
 \end{rem}

 \section{An application to the affine schemes of algebras.}\label{s:ap}

This section  shows that at least theoretically one may still be
able to answer the
 problem of Shafarevich discussed in the introduction by solving quadratic equations arising from
 the second order obstructions to deformations of algebras.
 In particular, we find that the vectors, contributing to the discrepancy
 between the tangent spaces to the scheme of associative multiplications and the
 scheme
 of the degree 3 nilpotent multiplications, can be easily ``killed'' by the
 obstructions in many cases.  Then Theorem~\ref{t:m} implies that  the tangent cones
 of the reduced schemes are the same, which means that the irreducible
 component of one variety is the component of the other one.

 Let us recall the notation from \cite{sh,am}.
 The affine scheme $C_n$ of all multiplications on a fixed $n$-dimensional vector
 space $V$
 over a field $K$ (${\rm char} K\ne2$),
 which represent associative commutative algebras,
 is given by the equations of commutativity and associativity:
 $$c_{ij}^k=c_{ji}^k,\qquad \sum_{s=1}^n c_{ij}^s c_{sk}^l=\sum_{s=1}^n c_{is}^l
 c_{jk}^s$$
 in the structure constants of multiplication
 $e_ie_j=\sum_{k=1}^n c_{ij}^k e_k$
 of the basis $\{e_1,\ldots,e_n\}$.
 Similar equations determine the affine scheme $A_n$ of commutative  nilpotent
 degree 3
 multiplications.  The reduced schemes associated to the above schemes are denoted
 $C_n^\red$ and $A_n^\red$, respectively.
 From \cite{sh}, the irreducible components of $A_n^\red$ have a very simple
 description
 $$A_{n,r}=\{N\in A_n^\red |\, \dim N^2\le r\le \dim {\rm Ann}_N N\},$$
 $1\le r\le(n-1)(n-r+1)/2$, where $N$ denotes the algebra represented by the
 corresponding
 multiplication, $N^2$ is the square of the algebra and $\rm Ann$ is the
 annihilator.

 Let $d:=n-r$, then, for $r=1,2$ and $r>(d^2-1)/3$, Shafarevich in \cite{sh} showed that
 $A_{n,r}$ is not a component of
 $C_n^\red$, by constructing
 a line, contained in  $C_n^\red$ but not in $A_n^\red$, through a point of
 $A_{n,r}$.
 For other $r$, one has to compare
 the tangent spaces to the non-reduced schemes $A_n$ and $C_n$.
 The smooth set of $A_n^\red$ is the union of
 $$U_{n,r}=\{N\in A_n^\red |\, \dim N^2= r= \dim {\rm Ann}_N N\}.$$
 If $W\subset V$ is a subspace of dimension $r$, then
 the space $S_{n,r}=L(S^2(V/W),W)$ of all linear maps from the symmetric product of
 $V/W$ to $W$
 is naturally included as an affine subspace into  $A_{n,r}$. The group $G={\rm
 GL}(V)$ acts
 on $C_n$ and $G S_{n,r}=A_{n,r}$.
 Then, the tangent space to $A_{n,r}$ at a point $N\in S_{n,r}$ is
 $$T_N A_{n,r}=L(S^2(N/N^2),N^2)+T_N G N,$$
 where the tangent space $T_N G N$ to the orbit is the space of
 maps from $L(S^2N,N)$ given by the coboundaries (in a Hochschild complex,
 see \cite{am}) $x\circ y= x\varphi(y)-\varphi(xy)+y\varphi(x)$
 for some linear map $\varphi:N@>>>N$.

 The tangent space $T_N C_n$
 always includes  $T_N A_{n,r}+F$, where
 the space $F$ consists of $x\circ y=xf(y)+yf(x)$ for some $f\in L(N/N^2,K)\subset
 L(N,K)$.
 To explain this difference
 Shafarevich embedded the scheme $C_n$ into the scheme $\widetilde{C}_n$
 of commutative associative multiplications on a $(n+1)$-dimensional space which
 represent algebras with a unit $e$:
 $$C_n\hookrightarrow\widetilde{C}_n, \qquad N\mapsto N\oplus Ke.$$
 In this situation, the subgroup $\widetilde G$ of $GL(V\oplus Ke)$ which fixes the
 unit
 acts on $\widetilde{C}_n$. It turns out that
 $T_N GN+F$ is the tangent space to
 $\widetilde{G} S_{n,r}$, whence  the equality
 \begin{equation}\label{e:tan}
 T_N \widetilde{C}_n=T_N A_{n,r}+F
 \end{equation}
 was enough to conclude that $A_{n,r}$ is the component of $C_n^\red$.
 The equality was shown in \cite{sh} for $3\le r\le (d+1)(d+2)/6$,
 and this also holds $r=(5d^2-8d)/16$ and $d$ divisible by 4
 by \cite[Section~2.1]{am}.

 The original Shafarevich's method is very special to this situation and
 is  impractical
 to use it in general to explain the nilpotents in the structure sheaf of a scheme.
 We expect that the equality (\ref{e:tan}) holds for all $3\le r<(d^2-1)/3$,
 which would leave only one case $r=(d^2-1)/3$ unsettled. However, in this last case
 the tangent space to $C_n$ at a generic point in $A_{n,r}$ is too big:
 \begin{equation}\label{e:for}
 T_N C_n= L(S^2(N/N^2),N^2)+T_N G N+L(S^2N_1,N_1),
 \end{equation}
 where $N=N_1\oplus N^2$ is a fixed decomposition.
 To show this one have to use Lemma~1 in \cite{am}:
 a  generic
 algebra $N\in A_{n,r}$, for $r\le(d^2-1)/3$, can be given by $d$ generators and a
 $(d(d+1)/2-r)$-dimensional space  of  homogeneous degree 2 relations among the
 generators,
 so that
 the nilpotence of the third degree follows from this relations.
 But, for $r=(d^2-1)/3$, this condition implies that there are no nontrivial
 relations among the
 homogeneous degree 2 relations, i.e., all of the relations $\sum_i n_iz_i=0$ in
 $S^3N_1$,
 $n_i\in N_1$, $z_i\in Z:=ker(\mu:S^2N_1@>>>N^2)$ ($\mu$ is the multiplication
 on $N$),  are induced by a trivial in $N_1\otimes Z$ element
 $\sum_i n_i\otimes z_i$. Indeed, $\dim S^3N_1=d(d+1)(d+2)/6$, while the number of
 equations
 that we get from the homogeneous degree 2 relations is equal to
 $$\dim N_1\otimes Z=d(d(d+1)/2-(d^2-1)/3)=d(d+1)(d+2)/6.$$
 So, if there is a nontrivial relation among the
 homogeneous degree 2 relations, there would be not enough equations to deduce
 $N^3=0$.
 Since the only restriction in \cite[Lemma~1]{am} for a vector from
 $L(S^2N_1,N_1)\subset L(S^2N,N)$ to be tangent to the scheme $C_n$
 was arising from the nontrivial relations,
 we conclude (\ref{e:for}).

 Anan'in suggested in \cite{am} to use the second order obstructions to deformations
 of algebras
 to eliminate the excessive space   $L(S^2N_1,N_1)$. To be precise, suppose that
 $N\in A_{n,r}$
 is a smooth point of $C_n^\red$ and $\circ\in L(S^2N_1,N_1)$. Then, if $\circ$ is
 tangent to $C_n^\red$, the sum $\circ+L(S^2(N/N^2),N^2)$ also lies in the tangent
 space
 $T_N C_n^\red$. A deformation of a commutative algebra on a vector space $V$
 with multiplication $x\cdot y$ can be thought as
 a smooth curve $$x\cdot y+(x\circ y)t+(x\star y)t^2+\cdots$$ in
 the affine space $L(S^2V,V)$, where $t$ is a local parameter.
 It is a well know fact in algebraic geometry that
 through a smooth point $p$ of a variety $X$ we can  find
 a smooth curve inside $X$, whose tangent vector at $p$ is a given one in $T_p X$.
 Therefore, the smooth curves in $C_n^\red$ with tangent vectors
 $\circ+L(S^2(N/N^2),N^2)$
   give a lot of equations (the second order obstructions) arising from the
 associativity of
 multiplication:
 $$(x\tilde\circ y)\tilde\circ z-x\tilde\circ(y\tilde\circ z)=
 x(y\star z)-(xy)\star z+x\star(yz)-(x\star y) z$$
 for some $\star\in L(S^2N,N)$, where $\tilde\circ\in\circ+L(S^2(N/N^2),N^2)$.
 A nice thing about these quadratic equations on $\circ$ is that one can linearize
 them:
 $$(x\circ y)* z+(x*y)\circ z-x\circ(y* z)-x*(y\circ z)=
 x(y\star z)-(xy)\star z+x\star(yz)-(x\star y) z$$
 for all $*\in L(S^2(N/N^2),N^2)$.
 These equations have been used to prove Theorem~1 in \cite{am} which implies that
 $A_{n,r}$ is the component of $C_n^\red$ for $d(d+1)/9\le r\le[d/3](d-3)$
 (almost all cases that are not covered by \cite{sh})
 if a certain algebra $N\in A_{n,r}$ is a smooth point
 of $C_n^\red$.
 Unfortunately, it seems impossible to prove that $N$ is the smooth point unless we
 know the tangent space at $N$.

 Now, we can apply Section~\ref{s:pc}.
 Let $N\in A_{n,r}$ with multiplication $x\cdot y$ and $\circ\in TC_N C_n^\red$.
 By Theorem~\ref{t:m}, there exists a smooth plane curve $C$
 $$x\cdot y+(x\circ y)t+(x\star y)t^2$$
 through $N$
 in the affine space $L(S^2V,V)$ of commutative
 multiplications on the vector space $V$, such that
 $\circ$ is the tangent vector to the curve at $N$ and the multiplicity of
 intersection $I(p,C,C_n^\red)\ge3$. Hence, the associativity equations in the
 structure
 constants imply
 \begin{equation}\label{e:ob}
 (x\circ y)\circ z-x\circ(y\circ z)=
 x(y\star z)-(xy)\star z+x\star(yz)-(x\star y) z
 \end{equation}
 for some $\star\in L(S^2N,N)$. Note that we do not assume the smoothness condition
 at $N$.

 We want to deduce  the effective part of the obstructions to $\circ$.
 As in \cite[\S1]{am}, decompose $\circ$ and $\star$ into the sum of linear maps
 $f_{ij}^k$ and $g_{ij}^k$ in $L(N_i\otimes N_j,N_k)$, respectively,
 where $1\le i,j,k\le2$ and $N=N_1\oplus N_2$,
 $N_2:=N^2$.
 For a sufficiently generic $N\in A_{n,r}$, besides commutativity,
 $f_{ij}^k$
  satisfy
 $$f_{12}^1=f_{22}^1=f_{22}^2=0,$$
 \begin{equation}\label{e:co}
 xf_{11}^1(y,z)+f_{12}^2(x,yz)=f_{12}^2(z,xy)+f_{11}^1(x,y)z,
 \end{equation}
 by \cite[\S1]{am}. Moreover,
 $f_{12}^2$ satisfying (\ref{e:co}) is unique for
 a given $f_{11}^1$, while its existence condition is described in
 \cite[Lemma~1]{am}.
 The missing $f_{11}^2\in L(S^2(N/N^2),N^2)$ is always tangent
 to $C_n$.

 From (\ref{e:ob}) we get
 \begin{equation}\label{e:ob1}
 f_{11}^1(f_{11}^1(x,y),z)-f_{11}^1(x,f_{11}^1(y,z))=
 -g_{12}^1(z,xy)+g_{12}^1(x,yz),
 \end{equation}
 \begin{equation}\label{e:ob2}
 f_{12}^2(x,f_{12}^2(y,zt))-f_{12}^2(y,f_{12}^2(x,zt))=
 xg_{12}^1(y,zt)-yg_{12}^1(x,zt)
 \end{equation}
 and
 $$f_{12}^2(f_{11}^1(x,y),zt)-f_{12}^2(x,f_{12}^2(y,zt))=xg_{12}^1(y,zt)-
 g_{22}^2(xy,zt),$$
 where $x,y,z,t\in N_1$.
 The last
 equation determines $g_{22}^2$, and, one can check that
 commutativity of $g_{22}^2$ (surprisingly)
 follows  from  (\ref{e:co}),
 (\ref{e:ob1}) and (\ref{e:ob2}).
 So, the only restrictions for $\circ$ to lie in the tangent cone to $C_n$
 are (\ref{e:co}),
 (\ref{e:ob1}) and (\ref{e:ob2}). The same argument as in \cite[\S1]{am}
 shows that $g_{12}^1$, satisfying (\ref{e:ob1}), is unique  for a given
 $f_{11}^1$, and the existence condition is similar to \cite[Lemma~1]{am}.
We have the following result:

\begin{thm}\label{t:m1}
The variety $A_{n,r}$ is an irreducible component
 of $C_n^\red$ when $r\le(d^2-1)/3$ if and only if,
 for some sufficiently generic $N\in A_{n,r}$,
a solution $\circ\in L(S^2N,N)$  to (\ref{e:co}),
 (\ref{e:ob1}) and (\ref{e:ob2}) has its part
$f^1_{11}=0$ on the kernel of the multiplication
$\mu:S^2N_1@>>>N^2$ of the algebra $N$.
\end{thm}

\begin{pf}
If $A_{n,r}$ is an irreducible component
 of $C_n^\red$, then, for some sufficiently generic $N\in A_{n,r}$,
 the tangent cone to $C_n^\red$ at $N$ coincides with
$$T_N A_{n,r}=L(S^2(N/N^2),N^2)+T_N G N.$$ But we know that the
solution to (\ref{e:co}),
 (\ref{e:ob1}) and (\ref{e:ob2}) gives rise to a vector $\circ$ from the
 tangent cone. Hence, $f^1_{11}=0$ on the kernel of the multiplication
$\mu:S^2N_1@>>>N^2$, because the vectors of $T_N G N$ are of the
form $x\circ y= x\varphi(y)-\varphi(xy)+y\varphi(x)$
 for some linear map $\varphi:N@>>>N$, while the vectors from
 $L(S^2(N/N^2),N^2)$ clearly satisfy the property.

 Conversely, suppose that a vector $\circ\in L(S^2N,N)$ belongs to the tangent
 cone to the scheme $C_n$ at $N$. So, by the discussion above, it
 satisfies
  (\ref{e:co}),
 (\ref{e:ob1}) and (\ref{e:ob2}). Hence,  $\circ$ has
$f^1_{11}=0$ on the kernel of the multiplication of some
sufficiently generic $N\in A_{n,r}$. Then $f_{11}^1$ with
$f_{12}^2(x,yz)=-xf_{11}^1(y,z)$ and
$g_{12}^1(x,yz)=-f_{11}^1(x,f_{11}^1(y,z))$ are unique solutions
to (\ref{e:co}), (\ref{e:ob1}) and (\ref{e:ob2}). This shows that
$\circ$ is actually from $T_N A_{n,r}$.
\end{pf}

 We will finish this section showing the result of
 \cite[Theorem~2]{am} with application of the above theorem and
 without the use of  Shafarevich's embedding $C_n\hookrightarrow\widetilde{C}_n$.

\begin{cor}
$A_{n,r}$ is a  component
 of $C_n^\red$ for the values $r=(5d^2-8d)/16$, $d$ is divisible by 4.
\end{cor}

\begin{pf}
 In \cite[Section~2.1]{am}, it was already shown that
 $T_N{C}_n=T_N A_{n,r}+F$ for a sufficiently generic algebra $N$.
 This implies that the part $f_{11}^1$ of a vector from $T_N{C}_n$ is determined by
 $f(x)y+f(y)x$, for some $f\in L(N/N^2,K)$, on the kernel of the multiplication
 $\mu$ on $N$.
 The algebra $N$ had the property that the square of all $d$ generators is zero and
 for each generator $u$ there was a distinct generator $v$
 such that their product $uv$ also vanishes in $N$.
 Taking $x=y=u$ and $z=v$ in (\ref{e:ob1}), we have
 $$f_{11}^1(2f(u)u,v)-f_{11}^1(u,f(u)v+f(v)u)=0.$$
 Since $uv$ and $uu$ are in the kernel of multiplication $\mu$, we further get
 $$2f(u)(f(u)v+f(v)u)-f(u)(f(u)v+f(v)u)-f(v)(2f(u)u)=0,$$
 whence $f(u)^2v-f(u)f(v)u=0$. Therefore, $f(u)=0$ for all generators, and
 $f_{11}^1=0$
 on the kernel of $\mu$ as it was required in Theorem~\ref{t:m1}.
\end{pf}

 \end{document}